\documentclass[11pt]{amsart}
\usepackage{amsmath,amsfonts,amsthm,amsopn,amssymb,mathtools,stmaryrd}
\usepackage{cite,marginnote}
\usepackage{pdfsync}

\usepackage{color,enumitem,graphicx}
\usepackage[colorlinks=true,urlcolor=blue,
citecolor=red,linkcolor=blue,linktocpage,pdfpagelabels,
bookmarksnumbered,bookmarksopen]{hyperref}
\usepackage[english]{babel}

\usepackage[left=2.75cm,right=2.75cm,top=2.85cm,bottom=2.85cm]{geometry}

\usepackage[hyperpageref]{backref}
\usepackage[]{xcolor}

\def\sideremark#1{\ifvmode\leavevmode\fi\vadjust{\vbox to0pt{\vss
 \hbox to 0pt{\hskip\hsize\hskip1em
 \vbox{\hsize2.1cm\tiny\raggedright\pretolerance10000
  \noindent #1\hfill}\hss}\vbox to15pt{\vfil}\vss}}}%

\numberwithin{equation}{section}

\makeindex
\newtheorem{theorem}{Theorem}[section]
\newtheorem{proposition}[theorem]{Proposition}
\newtheorem{lemma}[theorem]{Lemma}
\newtheorem{remark}[theorem]{Remark}
\newtheorem{example}[theorem]{Example}
\newtheorem{corollary}[theorem]{Corollary}
\newtheorem{definition}[theorem]{Definition}

\newcommand{\bt}{\begin{theorem}}
\newcommand{\et}{\end{theorem}}
\newcommand{\bl}{\begin{lemma}}
\newcommand{\el}{\end{lemma}}
\newcommand{\bd}{\begin{definition}}
\newcommand{\ed}{\end{definition}}
\newcommand{\bc}{\begin{corollary}}
\newcommand{\ec}{\end{corollary}}
\newcommand{\bp}{\begin{proof}}
\newcommand{\ep}{\end{proof}}
\newcommand{\bx}{\begin{example}}
\newcommand{\ex}{\end{example}}
\newcommand{\bi}{\begin{exercise}}
\newcommand{\ei}{\end{exercise}}
\newcommand{\bo}{\begin{proposition}}
\newcommand{\eo}{\end{proposition}}
\newcommand{\br}{\begin{remark}}
\newcommand{\er}{\end{remark}}
\newcommand{\be}{\begin{equation}}
\newcommand{\ee}{\end{equation}}
\newcommand{\ba}{\begin{align}}
\newcommand{\ea}{\end{align}}
\newcommand{\bn}{\begin{enumerate}}
\newcommand{\en}{\end{enumerate}}
\newcommand{\bg}{\begin{align*}}
\newcommand{\eg}{\end{align*}}
\newcommand{\bcs}{\begin{cases}}
\newcommand{\ecs}{\end{cases}}
\newcommand{\bean}{\begin{eqnarray*}}
\newcommand{\eean}{\end{eqnarray*}}

\newcommand{\s}{\section}

\newcommand{\lab}{\label}

\title[A Hartree-Fock system with mass subcritical 
growth]{Existence of normalized solutions of a Hartree-Fock system with mass subcritical growth}

\author[H.\ Jin]{Hua Jin}
\author[Y. Y.\ Chang]{Yanyun Chang}
\author[M. Squassina]{Marco Squassina}
\author[J. J.\ Zhang]{Jianjun Zhang}

\address[H.\ Jin]{\newline\indent School of Mathematics
\newline\indent
China University of Mining and Technology
\newline\indent
Xuzhou, 221116, China}
\email{\href{mailto:huajin@cumt.edu.cn}{huajin@cumt.edu.cn}}

\address[Y. Y.\ Chang]{\newline\indent School of Mathematics
\newline\indent
China University of Mining and Technology
\newline\indent
Xuzhou, 221116, China}
\email{\href{mailto:TS21080001A31@cumt.edu.cn}{TS21080001A31@cumt.edu.cn}}

\address[M. Squassina]{\newline\indent Dipartimento di Matematica e Fisica
\newline\indent
Universit\`a Cattolica del Sacro Cuore
\newline\indent
Via dei Musei 41, Brescia, Italy}\email{\href{marco.squassina@unicatt.it}{marco.squassina@unicatt.it}}

\address[J. J.\ Zhang]{\newline\indent College of Mathematica and Statistics
\newline\indent
Chongqing Jiaotong University
\newline\indent
Chongqing 400074, China}
\email{\href{mailto:zhangjianjun09@tsinghua.org.cn}{zhangjianjun09@tsinghua.org.cn}}

\thanks{(1) Corresponding author: \texttt{marco.squassina@unicatt.it}}
\thanks{(2) H. Jin and Y. Chang were supported by the Fundamental Research Funds for the Central Universities (2019XKQYMS90). M. Squassina is supported by Gruppo Nazionale per l'Analisi Ma\-te\-ma\-ti\-ca, la Probabilit\`a e le loro Applicazioni. J. J. Zhang was supported by NSFC (No.12371109, 11871123).}

\date{}

\begin{document}

\begin{abstract}
In this paper, we are concerned with normalized solutions in $H_{r}^{1}(\mathbb{R}^{3}) \times H_{r}^{1}(\mathbb{R}^{3})$ for Hartree-Fock type systems with the form
\be\lab{ Hartree-Fock}
\left\{
\begin{array}{ll}
-\Delta u +\alpha \phi _{u,v} u=\lambda _{1} u+\left | u \right | ^{2q-2} u+\beta \left | v \right | ^{q} \left | u \right | ^{q-2} u , \\
-\Delta v +\alpha \phi _{u,v} v=\lambda _{2} v+\left | v\right | ^{2q-2} v+\beta \left | u \right | ^{q} \left | v \right | ^{q-2} v , \\
\int_{\mathbb{R}^{3}}\left | u \right | ^{2} {\rm d}x=a_{1} , \quad \int_{\mathbb{R}^{3}}\left | v \right | ^{2} {\rm d}x=a_{2} , \nonumber\\
\end{array}
\right .
\ee
where $$
\phi_{u, v}\left(x\right):=\int_{\mathbb{R}^{3}} \frac{u^{2}(y)+v^{2}(y)}{|x-y|} {\rm d}y \in D^{1,2}\left(\mathbb{R}^{3}\right).
$$
\\
Here $\alpha,\beta>0, a_1,a_2>0$ and $1<q<\frac{5}{3}$. By seeking the constrained global minimizers of the corresponding functional, we prove that the existence of normalized solutions to the system above for any $a_1,a_2>0$ when $1<q<\frac{4}{3}$ and for $a_1,a_2>0$ small when $\frac{4}{3}\le q < \frac{3}{2}$. The nonexistence of normalized solutions is also considered for $\frac{3}{2}\le  q < \frac{5}{3}$. Also, the orbital stability of standing waves is obtained under local well-posedness assumptions of the evolution problem.
\end{abstract}

\maketitle

%
\s{Introduction}
\renewcommand{\theequation}{1.\arabic{equation}}
\subsection{Background}
The Hartree-Fock type system
\be\lab{Hartree-Fock2}
\left\{
\begin{array}{ll}
-\Delta _{x} \Psi _{1} +\alpha \left[\left | x \right |^{-1} \ast \left(\Psi _{1}^{2}+ \Psi _{2}^{2}\right) \right]\Psi _{1}=i\partial _{_{t} } \Psi _{1}+\left | \Psi _{1} \right |^{2q-2} \Psi _{1}+\beta \left | \Psi _{2} \right |^{q}\left | \Psi _{1} \right |  ^{q-2} \Psi _{1}, \\
-\Delta _{x} \Psi _{2} +\alpha \left[\left | x \right |^{-1} \ast \left(\Psi _{1}^{2}+ \Psi _{2}^{2}\right) \right]\Psi _{2}=i\partial _{_{t} } \Psi _{2}+\left | \Psi _{2} \right |^{2q-2} \Psi _{2}+\beta \left | \Psi _{1} \right |^{q}\left | \Psi _{2} \right |  ^{q-2} \Psi _{2},\\
\Psi _{j} =\Psi _{j}\left(x,t\right)\in \mathbb{C} , \left(x,t\right)\in\mathbb{R} ^{3} \times \mathbb{R},j=1,2,\\
\end{array}
\right.
\ee
has received a lot of attention in recent years. For instance, it appears in the basic quantum, chemistry model of the small number of electrons interacting with static nuclear, see \cite{Bris-5,Lieb-19,Lions-20} and the references therein for details. This system consists of two Schr\"odinger equations, in which there are Coulomb interaction terms. The constant $\beta \in \mathbb{R} $ describes the interspecies scattering lengths. When $\beta> 0$, it indicates interspecies attraction and $\beta< 0$ indicates interspecies repulsion.


Such problem was initially introduced by Hartree in \cite{Hartree-15} by employing a set of specialized test functions, without explicitly considering the Pauli exclusion principle. Subsequently, Fock in \cite{Fock-11} and Slater in \cite{Slater-22} addressed the Pauli exclusion principle by selecting a distinct class of test functions known as Slater determinants. By doing so, they derived a system of $N$-coupled nonlinear Schr\"odinger equations

\be\lab{IN1}
-\frac{\hbar^{2}}{2 m} \Delta \psi_{k}+V_{\mathrm{ext}} \psi_{k}+\left(\int_{\mathbb{R}^{3}}|x-y|^{-1}\sum_{j=1}^{N}\left|\psi_{j}\left(y\right)\right|^{2} {\rm d}y\right) \psi_{k}+\left(V_{\mathrm{ex}} \psi\right)_{k}=E_{k} \psi_{k},
\ee
where  $\psi_{k}: \mathbb{R}^{3} \rightarrow \mathbb{C}, k=1, \ldots, N, V_{\text {ext }} $ is a given external potential, and
$$
\left(V_{\text {ex }} \psi\right)_{k}:=-\sum_{j=1}^{N} \psi_{j} \int_{\mathbb{R}^{3}} \frac{\psi_{k}\left(y\right) \bar{\psi}_{j}\left(y\right)}{|x-y|} {\rm d}y
$$
is the $ k$-th  component of the crucial exchange term  and $ E_{k} $ is the $ k$-th eigenvalue. For more details about the Hartree-Fock method  we refer to \cite{Froehlich-12, Bokanowski-4, Mauser-21, Danielesen-7} and references therein.

In this paper, our main interest is focused on the case of $N=2$ and assume the external potential has the following form
$$
\left(V_{\mathrm{ex}} \psi\right)=-\left(\begin{array}{c}
\left|\psi_{1}\right|^{2 q-2} \psi_{1}+\beta\left|\psi_{1}\right|^{q-2}\left|\psi_{2}\right|^{q} \psi_{1} \\
\left|\psi_{2}\right|^{2 q-2} \psi_{2}+\beta\left|\psi_{1}\right|^{q}\left|\psi_{2}\right|^{q-2} \psi_{2}
\end{array}\right),
$$
which is consistent with the assumptions in \cite{Avenia-8}. It leads us to investigate the system (\ref{Hartree-Fock2}). Since we are mainly interested in the existence of standing wave solutions to (\ref{Hartree-Fock2}), namely, solutions having the form of
\be\lab{standing wave}
\Psi _{1} \left(x,t\right)=e^{-i\lambda _{1}t } u\left(x\right), \Psi _{2} \left(x,t\right)=e^{-i\lambda _{2}t } v\left(x\right), \lambda _{1} , \lambda _{2}\in \mathbb{R},
\ee
it suffices to consider the following coupled elliptic equations with nonlocal interaction
\be\lab{Hartree-Fock}
\left\{
\begin{array}{ll}
-\Delta u +\alpha \phi _{u,v} u=\lambda _{1} u+\left | u \right | ^{2q-2} u+\beta \left | v \right | ^{q} \left | u \right | ^{q-2} u \quad \text{in} \quad \mathbb{R} ^{3},\\
-\Delta v +\alpha \phi _{u,v} v=\lambda _{2} \,v+\left | v \right | ^{2q-2} v\,+\beta \left | u \right | ^{q} \left | v \right | ^{q-2} v  \quad \text{in} \quad \mathbb{R} ^{3},\\
\end{array}
\right.
\ee
where
$$
\phi_{u, v}\left(x\right):=\int_{\mathbb{R}^{3}}\frac{u^{2}\left(y\right)+v^{2}\left(y\right)}{|x-y|} {\rm d}y \in D^{1,2}\left(\mathbb{R}^{3}\right)
$$
is the unique solution in  $D^{1,2}\left(\mathbb{R}^{3}\right)$ of
$$
-\Delta \phi_{u, v}=4\pi\left(u^{2} +v^{2} \right) \quad \quad \text{in} \quad  \mathbb{R}^{3} .
$$
System (\ref{Hartree-Fock}) is called a Schr\"odinger-Poisson type system, see \cite{Faraj-10}.

In \cite{Avenia-8},  the authors  first studied the system (\ref{Hartree-Fock}), where $\lambda _{1} , \lambda _{2}\in \mathbb{R} $ are fixed parameter. They dealt with the functional
\begin{align*}
    \mathcal{I}\left ( u,v \right )=&\frac{1}{2}||\nabla u||_2^2+\frac{1}{2}||\nabla v||_2^2+\frac{\alpha}{4}\int_{\mathbb{R}^{3}}\left (u^2+v^2\right ) \phi_{u,v}{\rm d}x\\
    &-\frac{1}{2}\left( \lambda _{1}\left \| u \right \|_{2}  ^{2} +\lambda _{2} \left \| v \right \|_{2}  ^{2}\right)-\frac{1}{2q}\left (||u||^{2q}_{2q}+||v||_{2q}^{2q}\right ) -\frac{\beta}{q}\int_{\mathbb{R}^{3}}|u|^q|v|^q {\rm d}x
\end{align*}
and looked for its critical points in $H_{r}^{1}(\mathbb{R}^{3}) \times H_{r}^{1}(\mathbb{R}^{3})$. In that direction, mainly by variational methods, they
showed the existence of semitrivial and vectorial ground states solutions depending on the parameters involved. In addition, the authors in \cite{Chen-26} considered the least energy solutions of Hartree-Fock systems when the nonlinearities are subcritical. However, nothing can be said a priori on the $L^{2}$-norm of solutions.

In recent years, the study of normalized solutions has attracted considerable attentions, that is, the desired solutions have a priori prescribed $L^{2}$- norm. Let us introduce some related results about the Schr\"odinger-Poisson equations
$$
-\Delta u+\phi_{u} u-|u|^{p-2} u=\omega u \quad \text { in } \mathbb{R}^{3},
$$
where
$\phi_{u}\left(x\right)=\int_{\mathbb{R}^{3}} \frac{|u(y)|^{2}}{|x-y|} {\rm d}y$ satisfies $-\Delta \phi _{u} =4\pi u^{2} $.
In the last decades, the existence and stability of normalized solutions have been studied by many authors. We refer the reader to \cite{Bellazzini-2,Bellazzini-1,Jeanjean-16,Snchez-23, Bellazzini-3} and the references therein.
 The usual way in studying such problem is to looking for the constrained critical points of the functional
$$
\mathcal{J}\left(u\right)=\frac{1}{2} \int_{\mathbb{R}^{3}}|\nabla u|^{2} \mathrm{~d} x+\frac{1}{4} \int_{\mathbb{R}^{3}}  \int_{\mathbb{R}^{3}}\frac{|u\left(x\right)|^{2}|u\left(y\right)|^{2}}{|x-y|} \mathrm{d} x \mathrm{~d} y-\frac{1}{p} \int_{\mathbb{R}^{3}}|u|^{p} \mathrm{~d} x
$$
on the constraint
$$
S(c)=\left\{u \in H^{1}\left(\mathbb{R}^{3}\right): \int_{\mathbb{R}^{3}}|u|^{2} \mathrm{~d} x=c\right\}.
$$
In \cite{Snchez-23}, the authors proved the existence of minimizers when $p=\frac{8}{3}$, and $c\in\left(0,c_{0}\right)$ for a suitable $c_{0}>0$. When $p\in\left(2,3\right)$, it was shown in \cite{Bellazzini-3} that a minimizer exists if $c>0$ is small enough.
In \cite{Bellazzini-2}, J. Bellazzini and G. Siciliano obtained the existence and stability only for sufficiently large $L^{2}$-norm in case $3<p<\frac{10}{3}$, in case $p=\frac{8}{3}$ for sufficiently small charges. In \cite{Jeanjean-16}, L. Jeanjean and T. Luo gave a threshold value of $c_{1}>0$ for existence and nonexistence by a detailed study of the function $c\to m(c):=\inf _{u\in S_c} \mathcal{J}\left(u\right)$ in the range $p\in [3,\frac{10}{3}]$. Also, they gave a nonexistence result of normalized solutions when $p=3$ for all $c>0$ and when $p=\frac{10}{3}$ for $c>0$ is small enough. In addition, when $p\in\left(\frac{10}{3},6\right)$, $m\left(c\right)=-\infty$ for all $c>0$. In \cite{Bellazzini-1}, the authors considered the mass supercritical case $p\in\left(\frac{10}{3},6\right)$. By virtue of a mountain-pass argument developed on $S\left(c\right)$, they showed that for $c>0$ small enough, $\mathcal{J}$ admits a critical point constrained on $S\left(c\right)$ at a strictly positive energy level and it is orbitally unstable.

As for the existence of normalized solutions to nonlinear Schr\"odinger system
$$
\left\{\begin{array}{l}
-\Delta u+\lambda_{1} u=\mu_{1}|u|^{p-2} u+ \alpha|u|^{\alpha-2}|v|^{\beta} u \quad \text { in } \mathbb{R}^{N}, \\
-\Delta v+\lambda_{2} v=\mu_{2}|v|^{q-2} \,v+\, \beta|u|^{\alpha}|v|^{\beta-2}v \quad \text { in } \mathbb{R}^{N}, \\
\int_{\mathbb{R}^{N}} u^{2}=a^{2}, \quad \int_{\mathbb{R}^{N}} v^{2}=b^{2},
\end{array}\right.
$$
we refer to \cite{Gou-13, Deng-30, Bartsch-1-1, Bartsch-2, Bartsch-3, Gou-19, Bartsch-1, Zhong-25}) and point out that no nonlocal terms are involved.
In \cite{Wang-24}, J. Wang and W. Yang studied the coupled nonlinear Hartree equations with nonlocal interaction
$$
\left\{\begin{array}{c}
-\Delta u+V_{1}(x) u=\lambda_{1} u+\mu_{1}\left(\int_{\mathbb{R}^{N}} \frac{u^{2}(y)}{|x-y|^{2}} {\rm d}y\right) u+\beta\left(\int_{\mathbb{R}^{N}} \frac{v^{2}(y)}{|x-y|^{2}} {\rm d}y\right) u, \quad x \in \mathbb{R}^{N}, \\
-\Delta v+V_{2}(x) v=\lambda_{2} v+\mu_{2}\left(\int_{\mathbb{R}^{N}} \frac{v^{2}(y)}{|x-y|^{2}} {\rm d}y\right) v+\beta\left(\int_{\mathbb{R}^{N}} \frac{u^{2}(y)}{|x-y|^{2}} {\rm d}y\right) v, \quad x \in \mathbb{R}^{N}.
\end{array}\right.
$$
In addition to prove the existence and nonexistence of normalized solutions, they also obtained a precise description of the concentration behavior of solutions to the system under certain type trapping potentials by proving some delicate energy estimates. Due to the influence of nonlocal terms, we should emphasize that it is more difficult to estimate the energy and obtain the compactness of the $(PS)$ sequence, which also leads to less research on such problems.

\subsection{Main results}
Motivated by these recent works above, we consider the existence of solutions to (\ref{Hartree-Fock}) satisfying the conditions
\be\lab{condition}
\int_{\mathbb{R}^{3}}\left | u \right | ^{2} {\rm d}x=a_{1}>0 \quad \text{and} \quad \int_{\mathbb{R}^{3}}\left | v \right | ^{2} {\rm d}x=a_{2}>0.
\ee

Define
$$
S\left(a_{1},a_{2}\right):=\left \{  \left(u,v\right)\in H_{r}^{1}\left(\mathbb{R}^{3}\right) \times H_{r}^{1}\left(\mathbb{R}^{3}\right):\left \| u \right \| _{2}^{2} =a_{1},\left \| v \right \| _{2}^{2} =a_{2} \right \}
$$
and a solution $\left(u,v\right)\in H_{r}^{1}\left(\mathbb{R}^{3}\right) \times H_{r}^{1}\left(\mathbb{R}^{3}\right)$ of (\ref{Hartree-Fock})-(\ref{condition}) can be obtained by seeking a critical point of the functional
$$
I\left(u,v\right)=\frac{1}{2}||\nabla u||_2^2+\frac{1}{2}||\nabla v||_2^2+\frac{\alpha}{4}\int_{\mathbb{R}^{3}}\left (u^2+v^2\right ) \phi_{u,v}{\rm d}x-\frac{1}{2q}\left (||u||^{2q}_{2q}+||v||_{2q}^{2q}\right ) -\frac{\beta}{q}\int_{\mathbb{R}^{3}}|u|^q|v|^q {\rm d}x
$$
constrained on $S\left(a_{1},a_{2}\right)$. The parameters  $\lambda _{1} , \lambda _{2}\in \mathbb{R} $ are no longer fixed, but appear as Lagrange multipliers. $I$ is a functional of $C^{1}$-class and bounded from below when $1<q<\frac{5}{3}$. Let
\be\lab{minimization problem}
m\left(a_{1}, a_{2}\right):=\inf _{\left(u, v\right) \in S\left(a_{1}, a_{2}\right)} I\left(u, v\right).
\ee
In the present paper, by analyzing the compactness of the minimizing sequence of the related constraint problem, we obtain the existence of the normalized solutions of system (\ref{Hartree-Fock}). The orbital stability and some nonexistence results are also considered.

We state the main results as follows.
\bt\lab{mainresult}
Assume  $\alpha, \beta>0$ and $1<q<\frac{5}{3} $.
\begin{enumerate}
\item[(i)] When $1<q<\frac{4}{3} $, problem (\ref{Hartree-Fock})-(\ref{condition}) admits a normalized solution for any  $a_{1},  a_{2}>0 $.
\item[(ii)] When $\frac{4}{3}\le q < \frac{3}{2}$, problem (\ref{Hartree-Fock})-(\ref{condition}) admits a normalized solution for $a_{1}, a_{2}>0$ small.
\item[(iii)] When $\frac{3}{2}< q < \frac{5}{3}$, problem (\ref{Hartree-Fock})-(\ref{condition}) admits no normalized solution  for $a_{1}, a_{2}>0$ small.
\item[(iv)] When $q=\frac{3}{2}$, $1\le \alpha <8\pi$ and $0<\beta<\alpha$, problem (\ref{Hartree-Fock})-(\ref{condition}) admits no normalized solution for any $a_{1}, a_{2}>0$.
\end{enumerate}
\et
Next, we consider the orbital stability of the set of minimizers. 
\begin{definition}
Let
$$
G\left(a_{1}, a_{2}\right)=\left\{\left(u, v\right)\in S(a_{1}, a_{2}): I\left(u, v\right)=m\left(a_{1}, a_{2}\right)\right\}.
$$
$G\left(a_{1},a_{2}\right)$ is called orbitally stable, if for every $\varepsilon > 0$, there exists $\delta > 0$ so that if the initial datum $\left(\Psi _{1} \left (\cdot, 0 \right ),\Psi _{2} \left (\cdot, 0 \right )\right)$ in the system (\ref{Hartree-Fock2}) satisfies
$$
\inf _{\left ( u ,v \right )\in G\left ( a_{1} ,a_{2}\right )  }\left \| (\Psi _{1} \left (\cdot, 0 \right ) ,\Psi _{2} \left (\cdot, 0 \right ))-\left(u,v\right) \right \|_{H^{1}}< \delta ,
$$
there holds that
$$
\inf _{\left ( u,v \right )\in G\left ( a_{1} ,a_{2}\right )  }\left \| \left(\Psi _{1} \left (\cdot,t  \right ) ,\Psi _{2} \left (\cdot,t  \right )\right)-\left ( u,v\right ) \right \|_{H^{1}}< \varepsilon , \quad \forall  t>0 ,
$$
where $\Psi _{i}\left(\cdot,t\right)$ $i=1,2$ is the solution of (\ref{Hartree-Fock2}) with initial datum $\left(\Psi _{1} \left (\cdot, 0 \right ) ,\Psi _{2} \left (\cdot, 0 \right )\right)$.
\end{definition}
\bt\lab{stability}
Let $q\in \left ( 1,\frac{3}{2}  \right )$. Then the set $G\left(a_{1},a_{2}\right)$
is orbitally stable.
\et
\subsection{Main difficulties and ideas}
The main difficulty of the problem is the compactness of the minimizing sequence with respect to $m\left ( a_{1} ,a_{2}  \right )$. In order to overcome this difficulty, the method in \cite{Gou-13} is adopted. We consider the problem in $H_{r}^{1}(\mathbb{R}^{3}) \times H_{r}^{1}(\mathbb{R}^{3})$. By establishing a weak subadditive inequality, the strong convergence of the minimizing sequence is obtained. For the non-existence results, we mainly obtain it by a delicate estimate of the nonlocal term and applying the fact that any critical point of $I(u,v)$ on $S(a_{1},a_{2})$ satisfies the identity $Q(u,v)=0$, where $Q(u,v)$ is defined in (\ref{1536}).
In addition, through the scaling transformation  $u_{\theta}\left(x\right)=\theta ^{2} u\left(\theta x \right)$, $v_{\theta}\left(x\right)=\theta ^{2} v\left(\theta x \right)$,  compared with the case of a single Schr\"odinger-Poisson equation, a new similar $L^{2}$-critical index $q=\frac{4}{3} $ appears in our study. That is, when $1<q<\frac{4}{3} $, $m\left(a_{1}, a_{2}\right)<0$ for any  $a_{1}, a_{2}>0$, but when $\frac{4}{3}\le q < \frac{3}{2}$, $m\left(a_{1}, a_{2}\right)<0$ only for sufficiently small $a_{1}, a_{2}$, when $\frac{3}{2}< q < \frac{5}{3}$, $m\left(a_{1}, a_{2}\right)<0$ only for sufficiently large $a_{1}, a_{2}$.

\subsection{Notation}
\begin{itemize}
\item Denote the norm of $ L^{p}\left(\mathbb{R}^{3}\right)$  by
$$
\|u\|_{p}:=\left(\int_{\mathbb{R}^{3}}|u|^{p}{\rm d}x\right)^{\frac{1}{p}},\,\,1 \leq p<\infty.
$$
\item $ H^{1}\left(\mathbb{R}^{3}\right)$  is the usual Sobolev space endowed with the norm
$$
\|u\|_{H^{1}}:=\left(\int_{\mathbb{R}^{3 }  }^{} \left | \nabla u  \right |^{2} {\rm d}x+\int_{\mathbb{R}^{3 }  }^{} \left |  u  \right |^{2} {\rm d}x\right)^{\frac{1}{2}}
$$
and
$$
H_{r}^{1}\left(\mathbb{R}^{3}\right):=\left\{u \in H^{1}\left(\mathbb{R}^{3}\right): u\left(x\right)=u\left(|x|\right)\right\}.
$$
\item $D ^{1,2}\left(\mathbb{R}^{3} \right):=\left \{ u\in L^{2^{*}} (\mathbb{R}^{3} ):\nabla u\in L^{2}\left(\mathbb{R}^{3} \right)\right \}$ with the norm $\left ( \int_{\mathbb{R}^{3 }  }^{} \left | \nabla u  \right |^{2} {\rm d}x \right ) ^{\frac{1}{2} }$ and  $$D_{r}^{1,2}\left(\mathbb{R}^{3}\right):=\left \{u \in D ^{1,2}\left(\mathbb{R}^{3} \right):u\left(x\right)=u\left(|x|\right)\right\}.$$
\item Denote by $'\rightharpoonup ' $ and $ '\to'$ weak convergence and strong convergence, respectively. 
\item $C$ represents various positive constants which may be different from line to line. 
\item The symbol $ o_{n}(1)$ is used to denote a quantity that goes to zero as $n\to +  \infty$.
\end{itemize}
This paper is organized as follows.
In Section 2, some preliminaries are introduced. Particularly, some results in \cite{Avenia-8} are recalled that will be used to get compactness. We also give the variational setting for our problem.
Section 3 is devoted to the proof of Theorem \ref{mainresult}, which is about the existence and nonexistence of normalized solutions of (\ref{Hartree-Fock}).
In Section 4, the orbital stability of the set of minimizers is established.

\s{Preliminary results}
Firstly, let us observe that the $C^{1}$ functional $I\left(u,v\right)$ is well-defined in $H^{1}\left(\mathbb{R}^{3}\right) \times H^{1}\left(\mathbb{R}^{3}\right)$.
For $1<q<\frac{5}{3} $, thanks to the H\"older inequality, there is $p>1$ with $2<qp, qp'\le 6 , p':=\frac{p}{p-1} $, hence
$$
\int_{\mathbb{R}^{3}  }^{}\left | u \right |^{q} \left | v \right |^{q} {\rm d}x\le \left \| u \right \| _{qp}^{q}\left \| v \right \| _{qp'}^{q} < \infty \quad for\; u,v\in H^{1}\left(\mathbb{R}^{3}\right).
$$

We now give an upper bound estimate for the nonlocal term.
\bl
There exist constants $C_{1},C_{2},C_{3}>0$ independent of u and v, such that for all $u,v\in H^{1} \left ( \mathbb{R} ^{3}  \right ) $,
$$
\int_{\mathbb{R}^{3}}\left(u^2+v^2\right)\phi_{u,v}{\rm d}x \le C_{1}\left \| u \right \|  _{2}^{\frac{4}{3} }  \left \| u \right \| _{\frac{8}{3} } ^{\frac{8}{3} }+C_{2}\left \| v \right \|  _{2}^{\frac{4}{3} }  \left \| v \right \| _{\frac{8}{3} } ^{\frac{8}{3} }+C_{3}\left \| u \right \|  _{2}^{\frac{2}{3} }  \left \| u \right \| _{\frac{8}{3} } ^{\frac{4}{3} } \left \| v \right \| _{2} ^{\frac{2}{3}}\left \| v \right \| _{\frac{8}{3}} ^{\frac{4}{3}}.
$$
\bp
Since $\phi _{u}\left ( x \right ) :=\int_{\mathbb{R}^{3}} \frac{\left | u \left ( y \right ) \right | ^{2} }{|x-y|} {\rm d}y \in D^{1,2}\left(\mathbb{R}^{3}\right)$ solves the equation
\be\lab{u}
-\Delta \phi _{u} =4\pi u^{2}  \quad \text{in} \quad  \mathbb{R}^{3},
\ee
multiplying (\ref{u}) by $\phi _{u}\left ( x \right )$ and integrating, we obtain
$$
4\pi \int_{\mathbb{R}^{3} }^{} u^{2} \phi _{u}{\rm d}x =\int_{\mathbb{R}^{3} }^{}\left | \nabla \phi _{u} \right |^{2}{\rm d}x.
$$
Recall the following inequality
$$
\int_{\mathbb{R}^{3} }^{} u^{2} \phi _{u}{\rm d}x\le C\left \| u \right \|  _{2}^{\frac{4}{3} }  \left \| u \right \| _{\frac{8}{3} } ^{\frac{8}{3} },
$$
then we have
\begin{align*}
 \int_{\mathbb{R}^{3} }^{} v^{2} \phi _{u}{\rm d}x&\le \left ( \int_{\mathbb{R}^{3} }^{}\left | \phi _{u} \right |^{6}{\rm d}x \right )^{\frac{1}{6}} \left ( \int_{\mathbb{R}^{3} }^{}\left | v \right |^{\frac{12}{5} }{\rm d}x \right )^{\frac{5}{6}}\\
&\le   C\left ( \int_{\mathbb{R}^{3} }^{}\left | \nabla \phi _{u} \right |^{2}{\rm d}x \right ) ^{\frac{1}{2} }\left ( \int_{\mathbb{R}^{3} }^{}\left | v \right |^{\frac{12}{5} }{\rm d}x \right )^{\frac{5}{6}}\\
&=2\sqrt{\pi} C\left ( \int_{\mathbb{R}^{3} }^{} u^{2} \phi _{u}{\rm d}x \right )^{\frac{1}{2}}\left ( \int_{\mathbb{R}^{3} }^{}\left | v \right |^{\frac{12}{5} }{\rm d}x \right )^{\frac{5}{6}}\\
&\le \tilde{C} \left \| u \right \|  _{2}^{\frac{2}{3} }  \left \| u \right \| _{\frac{8}{3} } ^{\frac{4}{3} } \left \| v \right \| _{2} ^{\frac{2}{3}}\left \| v \right \| _{\frac{8}{3}} ^{\frac{4}{3}}.
\end{align*}
Thus,

\begin{align*}
 \int_{\mathbb{R}^{3}}\left(u^2+v^2\right)\phi_{u,v}{\rm d}x&=\int_{\mathbb{R}^{3} }^{} u^{2} \phi _{u}{\rm d}x+\int_{\mathbb{R}^{3} }^{} v^{2} \phi _{v}{\rm d}x+\int_{\mathbb{R}^{3} }^{} u^{2} \phi _{v}{\rm d}x+\int_{\mathbb{R}^{3} }^{} v^{2} \phi _{u}{\rm d}x\\
&\le C_{1}\left \| u \right \|  _{2}^{\frac{4}{3} }  \left \| u \right \| _{\frac{8}{3} } ^{\frac{8}{3} }+C_{2}\left \| v \right \|  _{2}^{\frac{4}{3} }  \left \| v \right \| _{\frac{8}{3} } ^{\frac{8}{3} }+C_{3}\left \| u \right \|  _{2}^{\frac{2}{3} }  \left \| u \right \| _{\frac{8}{3} } ^{\frac{4}{3} } \left \| v \right \| _{2} ^{\frac{2}{3}}\left \| v \right \| _{\frac{8}{3}} ^{\frac{4}{3}}.
\end{align*}
\ep
\el
Next, We begin to show that the following properties hold, which are important for proving the convergence of the minimizing sequence $\left ( u_{n},v_{n} \right )$ with respect to
$m(a_{1}, a_{2})$.
\bl\lab{PR11} (see\cite{Avenia-8} Lemma 3.2)
Let  $q \in(1,3) $ and  $\left\{\left(u_{n}, v_{n}\right)\right\} \subset H_{r}^{1}(\mathbb{R}^{3}) \times H_{r}^{1}(\mathbb{R}^{3})$ be such that $\left(u_{n}, v_{n}\right) \rightharpoonup (u, v) $ in $H_{r}^{1}(\mathbb{R}^{3}) \times H_{r}^{1}(\mathbb{R}^{3})$ as $ n \rightarrow+\infty $. We have, as $ n \rightarrow+\infty $,
\begin{align*}
\phi_{u_{n}, v_{n}} & \rightarrow \phi_{u, v} \text { in }
D_{r}^{1,2}\left ( \mathbb{R} ^{3}  \right ) , \\
\int_{\mathbb{R}^{3} }^{} \left(u_{n}^{2}+v_{n}^{2}\right) \phi_{u_{n}, v_{n}} {\rm d}x& \rightarrow \int_{\mathbb{R}^{3} }^{} \left(u^{2}+v^{2}\right) \phi_{u, v}{\rm d}x, \\
\int_{\mathbb{R}^{3} }^{} \left|u_{n}\right|^{q}\left|v_{n}\right|^{q}{\rm d}x & \rightarrow \int_{\mathbb{R}^{3} }^{} |u|^{q}|v|^{q}{\rm d}x.
\end{align*}
\el
As it is usual for elliptic equations, the
solutions of (\ref{Hartree-Fock}) satisfy a suitable identity called Pohozaev identity, which can be found in \cite[Lemma 3.1]{Avenia-8}. Benefiting from this Pohozaev identity, our nonexistence results are obtained.
\bl\lab{PR1}
If $\left(u,v\right)$ is a solution of (\ref{Hartree-Fock}), then it satisfies the Pohozaev identity
\begin{align*}
P_{\lambda _{1} ,\lambda _{2}} \left(u,v\right)&=\frac{1}{2}\left(||\nabla u||_2^2+||\nabla v||_2^2\right)-\frac{3}{2}\left(\lambda _{1} \left \| u \right \|_{2}  ^{2} +\lambda _{2} \left \| v \right \|_{2}  ^{2}\right) +\frac{5\alpha}{4}\int_{\mathbb{R}^{3}}\left(u^2+v^2\right)\phi_{u,v}{\rm d}x \\
&-\frac{3}{2q}\left(||u||^{2q}_{2q}+||v||_{2q}^{2q}\right)-\frac{3\beta}{q}\int_{\mathbb{R}^{3}}|u|^q|v|^q{\rm d}x\\
&=0.
\end{align*}
\el
\s{Proof of Theorem \ref{mainresult}}
\renewcommand{\theequation}{2.\arabic{equation}}
Before proving the main theorem, some Lemmas are in order. The next lemma shows that the functional $I\left(u,v\right)$ is bounded from below on $S(a_1,a_2)$ when $1<q<\frac{5}{3}$.
\bl\lab{coercive}
If $1<q<\frac{5}{3}$, then for every $a_1,a_2>0$, the functional $I\left(u,v\right)$ is bounded from below and coercive on $S\left(a_1,a_2\right)$.
\bp
The Gagliardo-Nirenberg inequality
\begin{align*}
    ||u||_p \leq C_{N,p}||\nabla u||_2^{\frac{N\left ( p-2 \right ) }{2p} }||u||_2^{1-\frac{N\left ( p-2 \right ) }{2p} } \quad for\; u\in H^{1}\left(\mathbb{R}^{N}\right),
\end{align*}
which holds for $2\leq p \leq 2^*$ when $N \geq 3$, implies for $\left ( u,v \right ) \in S\left ( a_{1},a_{2} \right )$,
$$\int_{\mathbb{R}^{3}} |u|^q|v|^q{\rm d}x\leq ||u||_{qp}^{q}||v||_{qp'}^{q}\leq C||\nabla  u||_2^{\frac{3(pq-2)}{2p}}||\nabla v||_2^{\frac{3(p'q-2)}{2p'}},$$
where $p>1$, $\frac{1}{p'} +  \frac{1}{p} =1$, $2\leq qp ,qp'\leq 6$.

So we obtain
\begin{align*}
    I\left(u,v\right)&=\frac{1}{2}||\nabla u||_2^2+\frac{1}{2}||\nabla v||_2^2+\frac{\alpha}{4}\int_{\mathbb{R}^{3}}\left(u^2+v^2\right)\phi_{u,v}{\rm d}x-\frac{1}{2q}\left(||u||^{2q}_{2q}+||v||_{2q}^{2q}\right)-\frac{\beta}{q}\int_{\mathbb{R}^{3}}|u|^q|v|^q{\rm d}x \\
    &\geq \frac{1}{2}||\nabla u||_2^2+\frac{1}{2}||\nabla v||_2^2-\frac{1}{2q}\left(||u||^{2q}_{2q}+||v||^{2q}_{2q}\right)-\frac{\beta}{q} \int_{\mathbb{R}^{3}}|u|^q|v|^q{\rm d}x \\
    &\geq \frac{1}{2}||\nabla u||^2_2+\frac{1}{2}||\nabla v||^2_2-\frac{1}{2q}\left(C_{q,a_1}||\nabla u||_2^{3\left(q-1\right)}+C_{q,a_2}||\nabla v||_2^{3\left(q-1\right)}\right)-\frac{\beta C}{q}||\nabla u||_2^{\frac{3(pq-2)}{2p}}||\nabla v||_2^{\frac{3(p'q-2)}{2p'}} \\
    &=\frac{1}{2}||\nabla u||_2^2 +\frac{1}{2}||\nabla v||^2_2-C_1||\nabla
    u||_2^{3\left(q-1\right)}-C_2||\nabla v||_2^{3\left(q-1\right)}-C_3||\nabla u||_2^{\frac{3(pq-2)}{2p}}||\nabla v||_2^{\frac{3(p'q-2)}{2p'}}.
\end{align*}

As $1<q<\frac{5}{3}$, it followes that $0<3\left(q-1\right)<2$, $0<\frac{3\left(pq-2\right)}{2p}+\frac{3\left(p'q-2\right)}{2p'}<2$,
which ensures the boundedness of $I\left ( u,v \right )$ from below and the coerciveness on $S\left(a_1,a_2\right)$.
\ep
\el
Hereafter, we use the same notation $m\left(a_1,a_2\right)$ for  $a_1,a_2\ge0$ with either $a_1>0$ or $a_2>0$, namely, one component of $\left(a_{1}, a_{2}\right)$ maybe zero.

In what follows, we collect some basic properties of $m\left(a_1,a_2\right)$.
\bl\lab{Infimum estimation}
\begin{enumerate}
\item[(1)] Let $1<q<\frac{4}{3}$, for any $a_1,a_2\ge0$ with either $a_1>0$ or $a_2>0$,
$$
-\infty <m\left(a_1,a_2\right)<0
$$.
\item[(2)] If $\frac{4}{3}\le q < \frac{3}{2}$, there exist $\rho _{1}, \rho _{2}>0$ such that
$
-\infty <m\left(a_1,a_2\right)<0
$ for all $a_{1}\in \left(0,\rho _{1}\right), a_{2}\in \left(0,\rho _{2}\right)$. If $\frac{3}{2}< q < \frac{5}{3}$, then there exist $\rho _{3}, \rho _{4}>0 $ such that
$
-\infty <m\left(a_1,a_2\right)<0
$ for all $a_{1}\in \left(\rho _{3},+\infty\right), a_{2}\in \left(\rho _{4},+\infty\right)$.
\item[(3)] $m\left(a_1,a_2\right)$ is continuous with respect to  $a_1, a_2$.
\item[(4)] For any $a_1\ge b_1\ge 0,a_2\ge b_2\ge 0$,
$$
m\left(a_1,a_2\right)\le m\left(b_1,b_2\right)+m\left(a_1-b_1,a_2-b_2\right).
$$
\end{enumerate}
\bp
(1) It follows from Lemma \ref{coercive} that $I\left(u,v\right)$ is coercive and in particular $m\left(a_1,a_2\right)>-\infty$. We define $u_s\left(x\right)=e^\frac{3s}{2} u\left(e^sx\right)$, $v_s\left(x\right)=e^\frac{3s}{2} v\left(e^sx\right)$, so that $\left \| u_s\left(x\right) \right \| _{2}^{2} =\left \| u\left(x\right) \right \| _{2}^{2} $, $\left \| v_s\left(x\right) \right \| _{2}^{2} =\left \| v\left(x\right) \right \| _{2}^{2} $, then we have the following scaling laws,
$$
\left \| \nabla u_s\left(x\right) \right \| _{2}^{2} =e^{2s} \left \| \nabla u\left(x\right) \right \| _{2}^{2} ,
$$
$$
\left \| u_s\left(x\right)\right \| _{2q}^{2q} =e^{3s\left(q-1\right)}\left \| u\left(x\right) \right \| _{2q}^{2q},
$$
$$
\int_{\mathbb{R}^{3}}^{}  \left | u_s\left(x\right) \right | ^q\left | v_s\left(x\right) \right | ^q{\rm d}x=e^{3s\left(q-1\right)}\int_{\mathbb{R}^{3}}^{}\left | u\left(x\right) \right | ^q\left | v\left(x\right) \right | ^q{\rm d}x  ,
$$
$$
\int_{\mathbb{R}^{3}}^{} \left[\left(u_s\left(x\right)\right)^2+\left(v_s\left(x\right)\right)^2\right]\phi _{u_s,v_s}\left(x\right){\rm d}x=e^s\int_{\mathbb{R}^{3}}^{} \left[u\left(x\right)^2+v\left(x\right)^2\right]\phi _{u,v}\left(x\right){\rm d}x.
$$
Thus
\begin{align*}
  I\left(u_{s} ,v_{s}\right)=&\frac{e^{2s}}{2} \left(\left \| \nabla u \right \| _{2}^{2} +\left \| \nabla v \right \| _{2}^{2}\right)+\frac{\alpha e^s}{4}\int_{\mathbb{R}^{3}}^{}\left(u^2+v^2\right)\phi  _{u,v}\left(x\right) {\rm d}x \\
 &-\frac{e^{3s\left(q-1\right)}}{2q}\left(\left \| u \right \| _{2q}^{2q}+\left \| v \right \| _{2q}^{2q} \right)-\frac{\beta e^{3s\left(q-1\right)}}{q} \int_{\mathbb{R}^{3}}^{}\left | u \right |  ^q \left | v \right |  ^q{\rm d}x .
\end{align*}
 We notice that $0<3s\left(q-1\right)<s  $ for $1<q<\frac{4}{3}$, thus, for $s \to -\infty$, we have $I\left(u_{s} ,v_{s}\right)\to 0^{-}$,
 which prove the first claim.\\
 \qquad (2) When $\frac{4}{3}\le q < \frac{3}{2}$,  we set $u_{\theta}\left(x\right)=\theta ^{\frac{1}{2}-\frac{3}{2} r } u\left(\frac{x}{\theta ^{r} } \right)$, $v_{\theta}\left(x\right)=\theta ^{\frac{1}{2}-\frac{3}{2} r } v\left(\frac{x}{\theta ^{r} } \right)$, so that $\left \| u_{\theta}\left(x\right) \right \| _{2}^{2} =\theta\left \| u\left(x\right) \right \| _{2}^{2} $, $\left \| v_{\theta}\left(x\right) \right \| _{2}^{2} =\theta\left \| v\left(x\right) \right \| _{2}^{2} $, then the following scaling laws can be obtained,
$$
\left \| \nabla u_{\theta}\left(x\right) \right \| _{2}^{2} =\theta ^{1-2r}  \left \| \nabla u\left(x\right) \right \| _{2}^{2} ,
$$
$$
\left \| u_{\theta}\left(x\right) \right \| _{2q}^{2q} =\theta ^{\left(1-3r\right)q+3r} \left \| u\left(x\right) \right \| _{2q}^{2q},
$$
$$
\int_{\mathbb{R}^{3}}^{}  \left | u_{\theta}\left(x\right) \right | ^q\left | v_{\theta}\left(x\right) \right | ^q{\rm d}x=\theta ^{\left(1-3r\right)q+3r} \int_{\mathbb{R}^{3}}^{}\left | u\left(x\right) \right | ^q\left | v\left(x\right) \right | ^q{\rm d}x  ,
$$
$$
\int_{\mathbb{R}^{3}}^{} \left[\left(u_{\theta}\left(x\right)\right)^2+\left(v_{\theta}\left(x\right)\right)^2\right]\phi _{u_{\theta},v_{\theta}}\left(x\right){\rm d}x=\theta^{2-r}\int_{\mathbb{R}^{3}}^{} \left[u\left(x\right)^2+v\left(x\right)^2\right]\phi _{u,v}\left(x\right){\rm d}x.
$$
Therefore
\begin{align*}
  I\left(u_{\theta} ,v_{\theta}\right)=&\frac{1}{2}\theta ^{1-2r} \left(\left \| \nabla u \right \| _{2}^{2} +\left \| \nabla v \right \| _{2}^{2}\right)+\frac{\alpha}{4}\theta^{2-r}\int_{\mathbb{R}^{3}}^{}\left(u^2+v^2\right)\phi  _{u,v}\left(x\right) {\rm d}x \\
 &-\frac{1}{2q}\theta ^{(1-3r)q+3r}\left(\left \| u \right \| _{2q}^{2q}+\left \| v \right \| _{2q}^{2q} \right)-\frac{\beta}{q} \theta ^{(1-3r)q+3r}\int_{\mathbb{R}^{3}}^{}\left | u \right |  ^q \left | v \right |  ^q{\rm d}x .
\end{align*}

Note that for $r=-1$ , we get
\begin{align*}
  I\left(u_{\theta} ,v_{\theta}\right)=&\frac{1}{2}\theta ^{3} \left(\left \| \nabla u \right \| _{2}^{2} +\left \| \nabla v \right \| _{2}^{2}\right)+\frac{\alpha }{4}\theta^{3}\int_{\mathbb{R}^{3}}^{}\left(u^2+v^2\right)\phi  _{u,v}\left(x\right) {\rm d}x \\
 &-\frac{1}{2q}\theta ^{4q-3}\left(\left \| u \right \| _{2q}^{2q}+\left \| v \right \| _{2q}^{2q} \right)-\frac{\beta}{q} \theta ^{4q-3}\int_{\mathbb{R}^{3}}^{}\left | u \right |  ^q \left | v \right |  ^q{\rm d}x .
\end{align*}
Since $4q-3<3$ for $\frac{4}{3}\le q<\frac{3}{2}$, there holds for $\theta \to 0$, $I\left(u_{\theta},v_{\theta}\right)\to 0^{-}$. Thus, there exist $\rho _{1}, \rho _{2}>0$ such that
$
-\infty <m\left(a_1,a_2\right)<0
$ for all $a_{1}\in \left(0,\rho _{1}\right), a_{2}\in \left(0,\rho _{2}\right)$. If $\frac{3}{2}<q<\frac{5}{3}$, we have  $4q-3>3$, then for $\theta \to +\infty$,  $I\left(u_{\theta} ,v_{\theta}\right)\to 0^{-}.$ Thus, there exist $\rho _{3}, \rho _{4}>0 $ such that
$
-\infty <m\left(a_1,a_2\right)<0
$ for all $a_{1}\in \left(\rho _{3},+\infty\right), a_{2}\in \left(\rho _{4},+\infty\right)$. The second claim is completed.
 \\
 \qquad (3) We assume $ \left(a_{1}^{n} , a_{2}^{n} \right)=\left(a_{1},a_{2}\right)+ o_{n} (1)$. From the definition of $m\left(a_{1}^{n} ,a_{2}^{n} \right)$, for any $ \varepsilon >0$, there exists $\left(u_n,v_n\right)\in S\left(a_{1}^{n} ,a_{2}^{n} \right)$ such that
 \be\lab{ll1}
 I\left(u_n,v_n\right)\le m\left(a_{1}^{n} ,a_{2}^{n} \right)+\varepsilon.
 \ee
 Setting $\bar{u} _{n} :=\frac{u_n}{\left \| u_n \right \|_2 }a _{1}^{\frac{1}{2} } $, $\bar{v} _{n} :=\frac{v_n}{\left \| v_n \right \|_2 }a _{2}^{\frac{1}{2} } $, we have that $\left(\bar{u} _{n},\bar{v} _{n}\right)\in S\left(a_1,a_2\right)$ and
 \be\lab{ll2}
 m\left(a_1,a_2\right)\le I\left(\bar{u} _{n},\bar{v} _{n}\right)=I\left(u_n,v_n\right)+ o_{n} (1).
 \ee
 Combining(\ref{ll1})and (\ref{ll2}) we obtain
 \be\lab{ll3}
  m\left(a_1,a_2\right)\le m\left(a_{1}^{n} ,a_{2}^{n} \right)+\varepsilon +o_{n}(1).
 \ee
  Similarly, from the definition of $ m\left(a_1,a_2\right) $, for any $ \varepsilon >0$, there exists $\left(u,v\right)\in S\left(a_1 ,a_2\right)$ such that
 \be\lab{ll4}
 I\left(u,v\right)\le m\left(a_1 ,a_2\right)+\varepsilon.
 \ee
 Let $\bar{u} :=\frac{u}{\left \| u\right \|_2 }\left(a_{1}^{n} \right)^{\frac{1}{2} } $, $\bar{v} :=\frac{v}{\left \| v \right \|_2 }\left(a_{2}^{n} \right)^{\frac{1}{2} } $, then $\left(\bar{u},\bar{v}\right)\in S\left(a_{1}^{n},a_{2}^{n}\right)$ and
 \be\lab{ll5}
 m\left(a_{1}^{n},a_{2}^{n}\right)\le I\left(\bar{u},\bar{v}\right)=I\left(u,v\right)+ o_{n} (1).
 \ee
 Combining(\ref{ll4})and  (\ref{ll5}) we deduce that
 \be\lab{hh}
 m\left(a_{1}^{n} ,a_{2}^{n} \right) \le  m\left(a_1,a_2 \right)+\varepsilon +o_{n}(1).
 \ee
 Therefore, since $\varepsilon >0$ is arbitrary, according to (\ref{ll3}) and (\ref{hh}) we deduce that
 $$
  m\left(a_{1}^{n} ,a_{2}^{n} \right) =  m\left(a_1,a_2 \right)+\varepsilon +o_{n}(1).
 $$
 The third claim is obtained.
 \\
 (4) By density of  $C_{0}^{\infty } \left(\mathbb{R^N} \right)$ into $H_{}^{1} \left(\mathbb{R^N} \right)$, for any $ \varepsilon >0$, there exist $\left(\bar{\xi} _{1} ,\bar{\xi} _{2}\right),\left(\hat{\xi} _{1} ,\hat{\xi} _{2}\right) \in C_{0}^{\infty } \left(\mathbb{R^N} \right)\times C_{0}^{\infty } \left(\mathbb{R^N} \right)$ with
 $\left \| \bar{\xi _i}  \right \| _{2}^{2} =b_i$,$\left \| \hat{\xi _i}  \right \| _{2}^{2} =a_i-b_i$ for $i=1,2$ such that
 \be\lab{1}
 I\left(\bar{\xi} _{1} ,\bar{\xi} _{2}\right) \le m\left(b_1,b_2\right)+\frac{\varepsilon }{2} ,
 \ee
 \be\lab{2}
 I\left(\hat{\xi} _{1} ,\hat{\xi} _{2}\right) \le m\left(a_1-b_1,a_2-b_2\right)+\frac{\varepsilon }{2} .
 \ee
 We may assume that
 $$
  \left ( supp\bar{\xi} _{1}\cup supp\bar{\xi} _{2} \right )  \cap \left ( supp\hat{\xi} _{1}\cup supp\hat{\xi} _{2} \right ) =\emptyset,
 $$
  and
 $$
 \int_{\mathbb{R}^{3} }^{} \left(\bar{\xi_{1}}^{2}+\bar{\xi_{2}}^{2}\right)\phi_{\hat{\xi} _{1},\hat{\xi} _{2}}{\rm d}x=\int_{\mathbb{R}^{3}}^{}\int_{\mathbb{R}^{3}}^{}\frac{\left(\bar{\xi_{1} } ^{2}\left ( x \right )+\bar{\xi_{2}}^{2}\left ( x \right )\right)\left(\hat{\xi_{1} } ^{2}\left ( y \right )+\hat{\xi_{2}}^{2}\left ( y \right )\right)}{\left | x-y \right | } <\frac{2\varepsilon}{\alpha},
 $$
then for $i=1,2$
 $$
 \left \| \bar{\xi} _{i}+\hat{\xi} _{i} \right \| _{2}^{2} =\left \| \bar{\xi} _{i}  \right \| _{2}^{2}+\left \| \hat{\xi} _{i}  \right \| _{2}^{2}=b_i+\left(a_i-b_i\right)=a_i.
 $$
 It follows that
 $m\left(a_1,a_2\right)\le I\left(\bar{\xi} _{1}+\hat{\xi} _{1} ,\bar{\xi} _{2}+\hat{\xi} _{2}\right) $. Set $\xi_i =\bar{\xi} _{i}+\hat{\xi} _{i}$ , we have $\left \| \xi _i \right \| _{2}^{2} =a_i$ for $i=1,2$.
and
 \begin{align*}
 I\left(\xi_{1},\xi_{2}\right)&=\frac{1}{2} \left(\left \| \nabla\xi_{1}  \right \| _{2}^{2} +\left \| \nabla\xi_{2}  \right \| _{2}^{2}\right)+\frac{\alpha }{4}\int_{\mathbb{R}^{3} }^{} \left(\xi_{1}^{2}+\xi_{2}^{2}\right)\phi_{\xi_{1},\xi_{2}} {\rm d}x \\
 &\ \ \ -\frac{1}{2q}\left(\left \| \xi_{1} \right \| _{2q}^{2q} +\left \| \xi_{2} \right \| _{2q}^{2q}\right) -\frac{\beta }{q} \int_{\mathbb{R}^{3} }^{} \left | \xi_{1 }\right |  ^{q}\left | \xi_{2} \right |  ^{q}{\rm d}x \\
 &=\frac{1}{2} \left(\left \| \nabla\bar{\xi}_{1}   \right \| _{2}^{2} +\left \| \nabla\bar{\xi_{2}} \right \| _{2}^{2}\right)+\frac{\alpha }{4}\int_{\mathbb{R}^{3} }^{} \left(\bar{\xi_{1}}^{2}+\bar{\xi_{2}}^{2}\right)\phi_{\bar{\xi_{1}},\bar{\xi_{2}} }{\rm d}x  \\
 &\ \ \ -\frac{1}{2q}\left(\left \| \bar{\xi_{1} }\right \| _{2q}^{2q} +\left \| \bar{\xi_2 }\right \| _{2q}^{2q}\right) -\frac{\beta }{q} \int_{\mathbb{R}^{3} }^{} \left | \bar{\xi_{1} }\right |  ^{q}\left | \bar{\xi_{2}} \right |  ^{q}{\rm d}x \\
&\ \ \ +\frac{1}{2} \left(\left \| \nabla\hat{\xi} _{1}   \right \| _{2}^{2} +\left \| \nabla\hat{\xi} _{2}\right \| _{2}^{2}\right)+\frac{\alpha }{4}\int_{\mathbb{R}^{3} }^{} \left(\hat{\xi} _{1}^{2}+\hat{\xi} _{2}^{2}\right)\phi_{\hat{\xi} _{1},\hat{\xi} _{2} }{\rm d}x\\
&\ \ \ -\frac{1}{2q}\left(\left \| \hat{\xi} _{1}\right \| _{2q}^{2q} +\left \| \hat{\xi} _{2}\right \| _{2q}^{2q}\right) -\frac{\beta }{q} \int_{\mathbb{R}^{3} }^{} \left | \hat{\xi} _{1}\right |  ^{q}\left | \hat{\xi} _{2}\right |  ^{q}{\rm d}x \\
&\ \ \ +\frac{\alpha }{4} \left ( \int_{\mathbb{R}^{3} }^{} \left(\bar{\xi_{1}}^{2}+\bar{\xi_{2}}^{2}\right)\phi_{\hat{\xi} _{1},\hat{\xi} _{2} }{\rm d}x +\int_{\mathbb{R}^{3} }^{} \left(\hat{\xi} _{1}^{2}+\hat{\xi} _{2}^{2}\right)\phi_{\bar{\xi_{1}},\bar{\xi_{2}} }{\rm d}x \right ) \\
&\le I\left(\bar{\xi_{1}},\bar{\xi_{2}}\right)+I\left(\hat{\xi} _{1},\hat{\xi} _{2}\right)+\varepsilon.
 \end{align*}
Combining  (\ref{1}) and (\ref{2}), we obtain
$$m\left(a_1,a_2\right)\le I\left(\xi_{1},\xi_{2}\right) \le m\left(b_1,b_2\right)+m\left(a_1-b_1,a_2-b_2\right)+2\varepsilon ,
$$
thus
$$
m\left(a_1,a_2\right)\le m\left(b_1,b_2\right)+m\left(a_1-b_1,a_2-b_2\right).
$$
This completes the proof of the Lemma.
\ep
\el

\br\lab{rr1}
Note that if we set  $u_s\left(x\right)=e^\frac{3s}{2} u\left(e^sx\right)$, $s\in\mathbb{R}$ then $$
\phi _{u_{s} } \left(x\right)=\int_{\mathbb{R}^{3} }^{} \frac{e^{3s}\left | u\left ( e^{s}y \right )  \right | ^{2}  }{\left | x-y \right | }dy=e^{s}\phi _{u}   \left ( e^{s} x \right )    .
$$
\er
To obtain our nonexistence results, we use the fact that any critical point of $I\left(u,v\right)$ restricted to $S\left ( a_{1},a_{2} \right )$ satisfies $Q\left(u,v\right)=0$, where
\be\lab{1536}
\begin{aligned}
  Q\left(u,v\right):=&||\nabla u||_2^2+||\nabla v||_2^2+\frac{\alpha}{4}\int_{\mathbb{R}^{3}}\left(u^2+v^2\right)\phi_{u,v}{\rm d}x-\frac{3\left(q-1\right)}{2q}\left(||u||^{2q}_{2q}+||v||_{2q}^{2q}\right)\\
  &-\frac{3\beta\left(q-1\right)}{q}\int_{\mathbb{R}^{3}}|u|^q|v|^q{\rm d}x.
\end{aligned}
\ee
Indeed, we have the following Lemmas.
\bl\lab{nonexistence}
If $\left(u_{0},v_{0}\right)$ is a critical point of $I\left(u,v\right)$ on $S\left(a_{1},a_{2}\right)$, then $Q\left(u_{0},v_{0}\right)=0$.
\bp
First, we denote
\begin{align*}
  I_{\lambda _{1} ,\lambda _{2}}\left(u,v\right)=&\frac{1}{2}||\nabla u||_2^2+\frac{1}{2}||\nabla v||_2^2-\frac{\lambda _{1} }{2} \left \| u \right \| _{2}^{2} -\frac{\lambda _{2} }{2} \left \| v \right \| _{2}^{2} +\frac{\alpha}{4}\int_{\mathbb{R}^{3}}\left(u^2+v^2\right)\phi_{u,v}{\rm d}x \\
&-\frac{1}{2q}(||u||^{2q}_{2q}+||v||_{2q}^{2q})-\frac{\beta}{q}\int_{\mathbb{R}^{3}}|u|^q|v|^q{\rm d}x,
\end{align*}
here $\lambda _{1} , \lambda _{2}\in R $  and $I_{\lambda _{1} ,\lambda _{2}}\left(u,v\right)$ is the energy functional corresponding to the equation (\ref{Hartree-Fock}).

Clearly,
$$
I_{\lambda _{1} ,\lambda _{2}}\left(u,v\right)=I\left(u,v\right)-\frac{\lambda _{1} }{2} \left \| u \right \| _{2}^{2} -\frac{\lambda _{2} }{2} \left \| v \right \| _{2}^{2},
$$
and simple calculations imply that
$$
Q\left(u,v\right)=\frac{3}{2}\left \langle I'_{\lambda _{1},\lambda _{2} }\left(u,v\right),\left(u,v\right)  \right \rangle -P_{\lambda _{1} ,\lambda _{2}} \left(u,v\right).
$$

Now from Lemma 3.1 of \cite{Avenia-8}, we know that $P_{\lambda _{1} ,\lambda _{2}} \left(u,v\right)=0$ is a Pohozaev identity for the Hartree-Fock equation (\ref{Hartree-Fock}). In particular, any critical point $\left(u,v\right)$ of $ I_{\lambda _{1} ,\lambda _{2}}\left(u,v\right) $ satisfies $ P_{\lambda _{1} ,\lambda _{2}} \left(u,v\right)=0$. On the other hand, since $\left(u_{0},v_{0}\right)$ is a critical point of $I\left(u,v\right)$ on $S\left(a_{1},a_{2}\right)$, there exists a Lagrange multiplier  $\lambda _{1} ,\lambda _{2}\in R $, such that $I'\left(u_{0},v_{0}\right)=\lambda _{1}\left ( u_{0}, 0 \right )+\lambda _{2}\left ( 0,v_{0} \right )$. Thus, for any $\left(\varphi _{1} ,\varphi _{2}\right)\in H^{1}\left({\mathbb{R}^{3}}\right)\times H^{1}\left({\mathbb{R}^{3}} \right),$ we have
$$
\left \langle I'_{\lambda _{1},\lambda _{2} }\left(u_{0},v_{0}\right),\left(\varphi _{1} ,\varphi _{2}\right)  \right \rangle =\left \langle I'\left(u_{0} ,v_{0} \right) -\lambda _{1}\left ( u_{0}, 0 \right ) -\lambda _{2}\left ( 0,v_{0} \right ) ,\left(\varphi _{1} ,\varphi _{2}\right)\right \rangle =0,
$$
which shows that $\left(u_{0},v_{0}\right)$ is also a critical point of $I_{\lambda _{1} ,\lambda _{2}}\left(u,v\right)$. Hence $P_{\lambda _{1} ,\lambda _{2}} \left(u_{0},v_{0}\right)=0$ and $$\left \langle I'_{\lambda _{1},\lambda _{2}  }(u_{0},v_{0}),\left(u_{0},v_{0}\right)  \right \rangle =0,$$ $Q\left(u_{0},v_{0}\right)=0$ follows then.
\ep
\el

 Now, a delicate estimate of the nonlocal term is given, which is available to control the functional $I\left(u,v\right)$ and $Q\left(u,v\right)$.
\bl\lab{control}
When $\frac{3}{2} \le q \le 2$, for any $\varepsilon>0$, there are  constants $C_{1}$, $C_{2}>0$ depending on $q$, $\varepsilon$, such that for any $ \left(u,v\right)\in S\left(a_{1},a_{2}\right)$ ,
$$
\int_{\mathbb{R}^{3}}\left(u^2+v^2\right)\phi_{u,v} {\rm d}x\ge -\frac{1}{8\pi\varepsilon^{2}} \left(||\nabla u||_2^2+||\nabla v||_2^2\right)+C_{1}  \frac{\left \| u \right \|_{2q}^{\frac{2q}{4-2q} }  }{\left \| \nabla u \right \|_{2}^{\frac{3\left(2q- 3\right)}{4-2q} } \left \| u \right \|_{2}^{\frac{2q-3}{4-2q} }   } + C_{2} \frac{\left \| v \right \|_{2q}^{\frac{2q}{4-2q} }  }{\left \| \nabla v \right \|_{2}^{\frac{3\left(2q-3\right)}{4-2q} } \left \| v \right \|_{2}^{\frac{2q-3}{4-2q} }   }.
$$
\el
\bp
When $\frac{3}{2} \le q \le 2$ , by interpolation, we have
\be\lab{12}
\left \| u \right \| _{2q}^{2q} \le \left \| u \right \| _{3}^{3\left(4-2q\right)}\left \| u \right \| _{4}^{4\left(2q-3\right)}.
\ee
Since the $\phi_{u, v}\left(x\right)\in D^{1,2}\left(\mathbb{R}^{3}\right)$ solves the equation
\be\lab{F2}
-\Delta \phi_{u, v}=4\pi\left(u^{2} +v^{2} \right) \quad \quad \text{in} \quad  \mathbb{R}^{3},
\ee
on one hand, multiplying (\ref{F2}) by $\phi_{u, v}\left(x\right)$ and integrating, we obtain
\be\lab{F3}
4\pi\int_{\mathbb{R}^{3}}\left(u^2+v^2\right)\phi_{u,v}{\rm d}x=\int_{\mathbb{R}^{3} }^{} \left | \nabla \phi_{u, v}\left(x\right) \right | ^{2} {\rm d}x.
\ee
On the other hand, multiplying (\ref{F2}) by $\left | u \right | +\left | v \right | $ and integrating, we get for any $\eta >0$ ,
\begin{equation}
\begin{aligned}\lab{F4}
 4\pi \eta\int_{\mathbb{R}^{3}}^{} \left(u^{2}+v^{2}\right)\left(\left | u \right |+\left | v \right | \right ){\rm d}x &=-\eta\int_{\mathbb{R}^{3}}^{} \Delta \phi_{u, v}\left(x\right)\left(\left | u \right |+\left | v \right | \right ){\rm d}x \\
   &=\eta \int_{\mathbb{R}^{3}}^{} \nabla \phi_{u, v}\left(x\right)\nabla (\left | u \right |+\left | v \right | ){\rm d}x.\\
\end{aligned}
\end{equation}
It follows from Young inequality that for any $\varepsilon>0$,
\begin{equation}
\begin{aligned}\lab{Y}
 4\pi \eta\int_{\mathbb{R}^{3}}^{} \left(u^{2}+v^{2}\right)\left(\left | u \right |+\left | v \right | \right ){\rm d}x \le \varepsilon\int_{\mathbb{R}^{3}}^{}\left | \nabla \phi_{u, v}\left(x\right) \right | ^{2}{\rm d}x +\frac{\eta ^{2}}{4\varepsilon}\int_{\mathbb{R}^{3}}^{} \left |\nabla\left(  u + v  \right) \right |^{2} {\rm d}x .
 \end{aligned}
\end{equation}
Thus, taking $\eta=1$ in (\ref{Y}), combining (\ref{F3}) and (\ref{Y}), we obtain
\be\lab{IN1}
4\pi\int_{\mathbb{R}^{3}}^{} \left(u^{2}+v^{2}\right)\left(\left | u \right |+\left | v \right | \right ){\rm d}x \le 4\pi\varepsilon\int_{\mathbb{R}^{3}}\left(u^2+v^2\right)\phi_{u,v}{\rm d}x+\frac{1}{4\varepsilon}\int_{\mathbb{R}^{3}}^{} \left |\nabla\left(  u + v  \right) \right |^{2} {\rm d}x .
\ee
Clearly, we observe that
\be\lab{e}
\int_{\mathbb{R}^{3}}^{} \left(u^{2}+v^{2}\right)\left(\left | u \right |+\left | v \right | \right ){\rm d}x \ge \int_{\mathbb{R}^{3}}^{} \left(\left | u \right |^{3}+\left | v \right |^{3}  \right){\rm d}x.
\ee
Then, from (\ref{IN1}) and  (\ref{e}),
\be\lab{cc}
 \left \| u \right \| _{3}^{3} +\left \| v \right \| _{3}^{3} \le \varepsilon \int_{\mathbb{R}^{3}}\left(u^2+v^2\right)\phi_{u,v}{\rm d}x+\frac{1}{16\pi\varepsilon}\int_{\mathbb{R}^{3}}^{} \left |\nabla  u + \nabla v \right |^{2} {\rm d}x ,
\ee
is obtained.
By (\ref{cc})
\begin{equation}
\begin{aligned}\lab{13}
\int_{\mathbb{R}^{3}}\left(u^2+v^2\right)\phi_{u,v}{\rm d}x &\ge \frac{1}{\varepsilon} \left ( \left \| u \right \| _{3}^{3} +\left \| v \right \| _{3}^{3} \right )-\frac{1}{16\pi\varepsilon^{2}}\int_{\mathbb{R}^{3}}^{} \left |\nabla  u + \nabla v \right |^{2} {\rm d}x\\
&\ge  \frac{1}{\varepsilon} \left ( \left \| u \right \| _{3}^{3} +\left \| v \right \| _{3}^{3} \right )-\frac{1}{8\pi\varepsilon^{2}}\left(||\nabla u||_2^2+||\nabla v||_2^2\right).
\end{aligned}
\end{equation}
Now, using Gagliardo-Nirenberg's inequality, there exists a constant $C_{q}>0$, such that
\be\lab{11}
\left \| u \right \| _{4}^{4(2q-3)} \le C_{q}\left \| \nabla u \right \| _{2}^{3(2q-3)} \left \| u \right \| _{2}^{2q-3}.
\ee
Taking (\ref{11}) into (\ref{12}), we obtain
\be
\left \| u \right \| _{2q}^{2q} \le C_{q}\left \| u \right \| _{3}^{3\left(4-2q\right)}\left \| \nabla u \right \| _{2}^{3\left(2q-3\right)} \left \| u \right \| _{2}^{2q-3}.
\ee
Thus,
\be\lab{18}
\left \| u \right \| _{3}^{3}\ge\frac{\widetilde{C_{q}}\left \| u \right \| _{2q}^{\frac{2q}{4-2q} }}{\left \| \nabla u \right \|_{2}^{\frac{3\left(2q-3\right)}{4-2q} }\left \| u \right \|_{2}^{\frac{2q-3}{4-2q} } }.
\ee
It follows from (\ref{18}) and  (\ref{13}) that
\begin{align*}
\int_{\mathbb{R}^{3}}\left(u^2+v^2\right)\phi_{u,v}{\rm d}x \ge &\frac{\widetilde{C_{q}}\left \| u \right \| _{2q}^{\frac{2q}{4-2q} }}{\varepsilon\left \| \nabla u \right \|_{2}^{\frac{3\left(2q-3\right)}{4-2q} }\left \| u \right \|_{2}^{\frac{2q-3}{4-2q} } }+\frac{\widetilde{C'_{q}}\left \| v \right \| _{2q}^{\frac{2q}{4-2q} }}{\varepsilon\left \| \nabla v \right \|_{2}^{\frac{3\left(2q-3\right)}{4-2q} }\left \| v \right \|_{2}^{\frac{2q-3}{4-2q} } }-\frac{1}{8\pi\varepsilon^{2}}\left(||\nabla u||_2^2+||\nabla v||_2^2\right)\\
&=\frac{C_{1}\left \| u \right \| _{2q}^{\frac{2q}{4-2q} }}{\left \| \nabla u \right \|_{2}^{\frac{3\left(2q-3\right)}{4-2q} }\left \| u \right \|_{2}^{\frac{2q-3}{4-2q} } }+\frac{C_{2}\left \| v \right \| _{2q}^{\frac{2q}{4-2q} }}{\left \| \nabla v \right \|_{2}^{\frac{3\left(2q-3\right)}{4-2q} }\left \| v \right \|_{2}^{\frac{2q-3}{4-2q} } }-\frac{1}{8\pi\varepsilon^{2}}\left(||\nabla u||_2^2+||\nabla v||_2^2\right)\\
\end{align*}
Then, the proof is completed.
\ep

The estimate on the nonlocal term leads to a lower bound on $Q\left(u,v\right)$.
\bl\lab{lower bound}
When $\frac{3}{2} < q< \frac{5}{3}$ and $ \alpha, \beta>0$, for any $\varepsilon>0$, there are constants $C_{3}\left ( \varepsilon, q, \alpha, \beta \right )$, $C_{4}\left ( \varepsilon, q, \alpha, \beta \right )>0$, such that for any $\left(u,v\right)\in S\left(a_{1},a_{2}\right)$,
\be\lab{223}
Q\left(u,v\right)\ge \frac{32\pi\varepsilon^{2}-\alpha }{32\pi\varepsilon^{2}} \left ( ||\nabla u||_2^2+||\nabla v||_2^2 \right ) -C_{3} \left \| \nabla u \right \|_{2}^{3} a_{1}^{\frac{1}{2} } -C_{4} \left \| \nabla v \right \|_{2}^{3} a_{2}^{\frac{1}{2} }.
\ee
\el
\bp
By Lemma \ref{control}, for any $\varepsilon>0$, there are constants $C_{1}>0$, $C_{2}>0$ depending on $\varepsilon, q$, such that, for any $ \left(u,v\right)\in S\left(a_{1},a_{2}\right)$, $\alpha, \beta>0$, there holds
\begin{equation}
\begin{aligned}\lab{224}
Q\left(u,v\right)\ge &\frac{32\pi\varepsilon^{2} -\alpha}{32\pi\varepsilon^{2}} \left(||\nabla u||_2^2+||\nabla v||_2^2\right)+\alpha C_{1} \frac{\left \| u \right \|_{2q}^{\frac{2q}{4-2q} }  }{\left \| \nabla u \right \|_{2}^{\frac{3\left(2q-3\right)}{4-2q} } \left \| u \right \|_{2}^{\frac{2q-3}{4-2q} }   } \\
&+ \alpha C_{2} \frac{\left \| v \right \|_{2q}^{\frac{2q}{4-2q} }  }{\left \| \nabla v \right \|_{2}^{\frac{3\left(2q-3\right)}{4-2q} } \left \| v \right \|_{2}^{\frac{2q-3}{4-2q} }   }-\frac{3\left(q-1\right)\left ( \beta +1 \right )}{2q}\left(||u||^{2q}_{2q}+||v||_{2q}^{2q}\right).
\end{aligned}
\end{equation}
To obtain (\ref{223}) from (\ref{224}), we introduce the auxiliary function
$$
f_{k_{1},k_{2}} \left(x_{1},x_{2}\right)=\frac{32\pi\varepsilon^{2} -\alpha}{32\pi\varepsilon^{2}}\left(k_{1}+k_{2}\right)+\alpha D_{1}x_{1}^{\frac{1}{4-2q}}+\alpha D_{2}x_{2}^{\frac{1}{4-2q}}-\frac{3\left(q-1\right)\left ( \beta +1 \right )}{2q}\left(x_{1}+x_{2}\right),\quad x_{1},x_{2}>0
$$
with $D_{1}=C_{1} \left(k_{1}^{\frac{3(2q-3)}{2\left(4-2q\right)}} \cdot a_{1}^{\frac{2q-3}{2\left(4-2q\right)}}\right)^{-1}$, and $D_{2}=C_{2} \left(k_{2}^{\frac{3(2q-3)}{2\left(4-2q\right)}} \cdot a_{2}^{\frac{2q-3}{2\left(4-2q\right)}}\right)^{-1}$.
The study of the auxiliary function will provide us with an estimate independent of $||u||^{2q}_{2q}$, $||v||^{2q}_{2q}$. Clearly,
$$
\frac{\partial f_{k_{1},k_{2}} \left(x_{1},x_{2}\right) }{\partial x_{i}} =\frac{\alpha}{4-2q} \cdot D_{i} \cdot x_{i}^{\frac{2q-3}{4-2q}}-\frac{3\left(q-1\right)\left ( \beta +1 \right )}{2q},
$$
$$
 \frac{\partial^2}{\partial x_i^{2}}f_{k_{1},k_{2}} \left(x_{1},x_{2}\right)=\frac{\alpha}{4-2q} \cdot \frac{2q-3}{4-2q} \cdot D_{i} \cdot x_{i}^{\frac{4q-7}{4-2q}}>0,\quad for\;  x_{1}, x_{2}>0, i=1,2.
$$
For convenience, we set $M:=\frac{3\left(q-1\right)\left ( \beta +1 \right )\left(4-2q\right)}{2q}$.
Therefore, $f_{k_{1},k_{2}} \left(x_{1},x_{2}\right)$ has the unique global minimum at
$$
\left[\bar{x_{1}},\bar{x_{2}}\right]=\left[\left(\frac{M}{\alpha D_{1}}\right)^{\frac{4-2q}{2q-3}},\left(\frac{M}{ \alpha D_{2}}\right)^{\frac{4-2q}{2q-3}}\right],
$$
and
\begin{align*}
f_{k_{1},k_{2}} \left(\bar{x_{1}},\bar{x_{2}}\right)&=\frac{32\pi\varepsilon^{2} -\alpha}{32\pi\varepsilon^{2}}\left(k_{1}+k_{2}\right)+\alpha D_{1}\left(\frac{M}{\alpha D_{1}}\right)^{\frac{1}{2q-3}}+\alpha D_{2}\left(\frac{M}{\alpha D_{2}}\right)^{\frac{1}{2q-3}}\\
&\ \ \ -\frac{3\left(q-1\right)\left ( \beta +1 \right )}{2q}\left[\left(\frac{M}{\alpha D_{1}}\right)^{\frac{4-2q}{2q-3}}+\left(\frac{M}{\alpha D_{2}}\right)^{\frac{4-2q}{2q-3}}\right]\\
&=\frac{32\pi\varepsilon^{2} -\alpha}{32\pi\varepsilon^{2}}\left(k_{1}+k_{2}\right)+\left (\alpha D_{1}\right )^{\frac{2q-4}{2q-3}} \cdot M^{\frac{1}{2q-3}}\cdot\left ( 1-\frac{1}{4-2q}  \right )\\
&\ \ \ +\left (\alpha D_{2}\right )^{\frac{2q-4}{2q-3}}\cdot M^{\frac{1}{2q-3}}\cdot\left ( 1-\frac{1}{4-2q}  \right )\\
&= \frac{32\pi\varepsilon^{2} -\alpha}{32\pi\varepsilon^{2}}\left(k_{1}+k_{2}\right)-\alpha^{\frac{2q-4}{2q-3}}\widetilde{C_{1}}\cdot M^{\frac{1}{2q-3} } \cdot \frac{2q-3}{4-2q}\cdot k_{1}^{\frac{3}{2}} \cdot a_{1}^{\frac{1}{2}}\\
&\ \ \ -\alpha^{\frac{2q-4}{2q-3}}\widetilde{C_{2}}\cdot M^{\frac{1}{2q-3} } \cdot \frac{2q-3}{4-2q}\cdot k_{2}^{\frac{3}{2}} \cdot a_{2}^{\frac{1}{2}}\\
&=\frac{32\pi\varepsilon^{2}-\alpha }{32\pi\varepsilon^{2}} \left ( ||\nabla u||_2^2+||\nabla v||_2^2 \right ) -C_{3} \left \| \nabla u \right \|_{2}^{3} a_{1}^{\frac{1}{2} } -C_{4} \left \| \nabla v \right \|_{2}^{3} a_{2}^{\frac{1}{2} }.
\end{align*}

Because of $f_{k_{1},k_{2}} \left(x_{1},x_{2}\right)\ge f_{k_{1},k_{2}} \left(\bar{x_{1}},\bar{x_{2}}\right)$ for all $x_{1}, x_{2}>0$, we get (\ref{223}).
\ep
Next, we are ready to give the proof of Theorem \ref{mainresult}.\\
$ Proof\ of \ Theorem$ \ref{mainresult}. Assume that $\left ( u_{n},v_{n} \right ) $ is a minimizing sequence with respect to
$m(a_{1}, a_{2})$, then $I\left ( u_{n},v_{n}   \right ) = m\left ( a_{1},a_{2}   \right )+ o _{n} \left ( 1 \right )$. By the coerciveness of $I\left ( u,v \right ) $ on $S\left ( a_{1},a_{2} \right ) $, the sequence $\left ( u_{n},v_{n} \right ) $ is bounded, and so, $\left ( u_{n} ,v_{n}  \right ) \rightharpoonup \left ( u,v \right ) $ in   $H_{r}^{1}(\mathbb{R}^{3}) \times H_{r}^{1}(\mathbb{R}^{3})$.
By the compactness of the embedding $H_{r}^{1}\left(\mathbb{R}^{N}\right)\subset L^{p}\left(\mathbb{R}^{N}\right)$ for $2<p<6$, Lemma \ref{PR11}, and the weak convergence, the following formulas hold
\begin{align*}
 u_{n} &\rightarrow  u \quad \text{in} \quad L^{2q}\left(\mathbb{R}^{3}\right),\\
 v_{n} &\rightarrow  v \quad \text{in} \quad L^{2q}\left(\mathbb{R}^{3}\right),\\
    \int_{\mathbb{R}^{3} }^{} \left(u_{n}^{2} +v_{n}^{2} \right) \phi_{u_{n}, v_{n}} {\rm d}x &\rightarrow \int_{\mathbb{R}^{3} }^{} \left (u^{2} + v ^{2} \right )  \phi_{u, v}{\rm d}x, \\
\int_{\mathbb{R}^{3} }^{} \left|u_{n}\right|^{q}\left|v_{n}\right|^{q}{\rm d}x & \rightarrow \int_{\mathbb{R}^{3} }^{} |u|^{q}|v|^{q}{\rm d}x,
\end{align*}
thus we have
\be\lab{2.7}
m\left(a_{1}, a_{2}\right)=\lim_{n \to \infty} I\left(u_{n}, v_{n}\right) \ge I\left(u,v\right).
\ee
Assume that $\left ( u_{n} ,v_{n}  \right ) \rightharpoonup \left ( u,v \right )=\left ( 0,0 \right ) $ in   $H_{r}^{1}(\mathbb{R}^{3}) \times H_{r}^{1}(\mathbb{R}^{3})$, it follows that $m\left(a_{1}, a_{2}\right)\ge 0$ which contradicts with $m\left(a_{1}, a_{2}\right)< 0$.
Note that if $\left\|u\right\|_{2}^{2}=a_{1}$ and $\left\|v\right\|_{2}^{2}=a_{2}$ we are done. Indeed, from the definition of $m\left(a_{1}, a_{2}\right)$, we deduce
$I\left(u, v\right) \ge m\left(a_{1}, a_{2}\right)$ this moment, this together with (\ref{2.7}) leads to
\be\lab{1507}
m\left(a_{1}, a_{2}\right)=I\left(u,v\right).
\ee
Therefore, combined with $I\left ( u_{n},v_{n}   \right )= m\left ( a_{1},a_{2}   \right )+ o _{n} \left ( 1 \right ) $, the strong convergence of $\left(u_{n},v_{n}\right)$ in $H_{r}^{1}\left(\mathbb{R}^{N}\right) \times H_{r}^{1}\left(\mathbb{R}^{N}\right)$ then directly follows.
Otherwise, we assume by contradiction that  $\left\|u\right\|_{2}^{2}:=b_{1}<a_{1}$ or $\left\|v\right\|_{2}^{2}:=b_{2}<a_{2}$. By definition $I\left(u,v\right)\ge m\left(b_{1},b_{2}\right)$ and thus it results from (\ref{2.7}) that
\be\lab{2.8}
m\left(b_{1},b_{2}\right) \le m\left(a_{1}, a_{2}\right).
\ee
At this point, by Lemma  \ref{Infimum estimation}, in  case $1<q<\frac{4}{3}$,
$
m\left(a_1-b_1,a_2-b_2\right)<0
$. In  case $\frac{4}{3}\le q < \frac{3}{2}$, then there are $\rho _{1},\rho _{2}>0 $ such that $m\left(a_1-b_1,a_2-b_2\right)<0$ for all $a_{1}\in \left(0,\rho _{1}\right),a_{2}\in \left(0,\rho _{2}\right)$. So we get
$$
m\left(a_1,a_2\right) > m\left(b_1,b_2\right)+m\left(a_1-b_1,a_2-b_2\right),
$$
which is a contradiction to Lemma\ref{Infimum estimation}(4) and Theorem \ref{mainresult} (1),(2) is proved.

Since there is $\left ( u,v \right ) \in S\left ( a_{1},a_{2}   \right )$ with $m\left(a_{1}, a_{2}\right)=I\left(u,v\right)$. By the Lagrange multiplier, there exist $\lambda _{1} ,\lambda _{2} \in \mathbb{R}$ such that
$$
I'\left ( u,v \right ) =\lambda _{1} \left ( u,0 \right ) +\lambda _{2} \left ( 0,v \right ).
$$
Therefore, we obtain the normalized solution $\left ( \lambda _{1},\lambda _{2},u,v  \right )$ of (\ref{Hartree-Fock})-(\ref{condition}) in $H_{r}^{1}(\mathbb{R}^{3}) \times H_{r}^{1}(\mathbb{R}^{3})$ for the above several cases.

We consider the non-existence for $\frac{3}{2}< q < \frac{5}{3}$. By contradiction, assuming that there are sequence $a_{1}^{n}\subset R^{+}$, $a_{2}^{n}\subset R^{+}$, with $a_{1}^{n}\to 0$, $a_{2}^{n}\to 0$, as $n\to \infty$, and $\left \{ \left ( u_{n},v_{n} \right )  \right \} \subset S\left(a_{1}^{n}, a_{2}^{n}\right)$ such that $\left( u_{n},v_{n}\right)\subset S\left(a_{1}^{n}, a_{2}^{n}\right)$ is a critical point of $I\left(u,v\right)$ restricted to $S\left(a_{1}^{n}, a_{2}^{n}\right)$ . Then, on the one hand, from Lemma \ref{nonexistence},
\begin{align*}
Q\left(u_{n},v_{n}\right)&=||\nabla u_{n}||_{2}^{2}+||\nabla v_{n}||_{2}^{2}+\frac{\alpha}{4}\int_{\mathbb{R}^{3}}\left(u_{n}^{2}+v_{n}^{2}\right)\phi_{u_{n},v_{n}}{\rm d}x\\
&-\frac{3\left(q-1\right)}{2q}\left(||u_{n}||^{2q}_{2q}+||v_{n}||_{2q}^{2q}\right)-\frac{3\beta\left(q-1\right)}{q}\int_{\mathbb{R}^{3}}|u_{n}|^{q}|v_{n}|^{q}{\rm d}x\\
&=0.
\end{align*}
Since $\alpha>0$, $\beta>0$ and $\frac{3}{2}<q<\frac{5}{3}$, naturally, we deduce
\be
||\nabla u_{n}||_{2}^{2}+||\nabla v_{n}||_{2}^{2}\le \frac{3\left(q-1\right)\left ( \beta +1 \right )}{2q}\left(||u_{n}||^{2q}_{2q}+||v_{n}||_{2q}^{2q}\right).
\ee
We have, from Gagliardo-Nirenberg's inequality, that for some $C_{1}>0$, and $C_{2}>0$,
\be
||\nabla u_{n}||_{2}^{2}+||\nabla v_{n}||_{2}^{2}\le C_{1}(a_{1}^{n})^{\frac{3-q}{2}}\left \| \nabla u_{n}  \right \| _{2}^{3\left(q-1\right) } +C_{2}(a_{2}^{n})^{\frac{3-q}{2}}\left \| \nabla v_{n}  \right \| _{2}^{3\left(q-1\right) }.
\ee
Because of $3\left(q-1\right)<2$, we obtain that
\be\lab{b}
||\nabla u_{n}||_{2}\to 0 \quad and \quad ||\nabla v_{n}||_{2}\to 0 \qquad as \qquad n \to \infty.
\ee
On the other hand, by Lemma \ref{lower bound}, it follows that there are constants $C_{3}\left ( \varepsilon, q, \alpha, \beta \right )$, $C_{4}\left ( \varepsilon, q, \alpha, \beta \right )>0$ such that
\be\lab{a}
\frac{32\pi\varepsilon^{2}-\alpha }{32\pi\varepsilon^{2}} \left ( ||\nabla u_{n}||_2^2+||\nabla v_{n}||_2^2 \right ) \le C_{3} \left \| \nabla u_{n} \right \|_{2}^{3} \left(a_{1}^{n}\right)^{\frac{1}{2} } +C_{4} \left \| \nabla v_{n} \right \|_{2}^{3} \left(a_{2}^{n}\right)^{\frac{1}{2} }.
\ee
According to the arbitrariness of $\varepsilon$, we can take $\varepsilon>\sqrt{\frac{\alpha }{32\pi} }$, then (\ref{a}) implies that
$$
||\nabla u_{n}||_{2}\to \infty \quad or \quad ||\nabla v_{n}||_{2}\to \infty \qquad as \qquad n \to \infty,
$$
which are contradictory to  (\ref{b}). Thus, we finish the proof of Theorem \ref{mainresult}(3).

Now when $q=\frac{3}{2}$, it is enough to prove that, for any $a_{1}.a_{2}>0$, there holds $Q\left(u,v\right)>0$ for all $\left(u,v\right)\in S(a_{1},a_{2})$.
Indeed, if $Q\left(u,v\right)>0$ holds true, we can conclude the nonexistence of minimizers directly from Lemma \ref{nonexistence}.

To check  $Q\left(u,v\right)>0$ for all $\left(u,v\right)\in S\left(a_{1},a_{2}\right)$, let $\eta=2$ in (\ref{F4}) and $\varepsilon=1$ in (\ref{Y}), then from (\ref{F3}) and (\ref{Y}), we get
\be\lab{F5}
\frac{\alpha}{4}\int_{\mathbb{R}^{3}}\left(u^2+v^2\right)\phi_{u,v}{\rm d}x\ge \frac{\alpha }{2} \int_{\mathbb{R}^{3}}\left(u^2+v^2\right)\left(\left | u \right |+\left | v \right | \right){\rm d}x-\frac{\alpha }{16\pi}\int_{\mathbb{R}^{3}}\left(\nabla (\left | u \right |+\left | v \right | )^{2} {\rm d}x\right) .
\ee
Thus, for any $\left(u,v\right)\in S\left(a_{1},a_{2}\right)$, $q=\frac{3}{2}$,
\begin{align*}
 Q\left(u,v\right)&=||\nabla u||_2^2+||\nabla v||_2^2+\frac{\alpha}{4}\int_{\mathbb{R}^{3}}\left(u^2+v^2\right)\phi_{u,v}{\rm d}x-\frac{1}{2}\left(||u||^{3}_{3}+||v||_{3}^{3}\right)-\beta\int_{\mathbb{R}^{3}}|u|^{\frac{3}{2}}|v|^{\frac{3}{2}}{\rm d}x  \\
 &\ge ||\nabla u||_2^2+||\nabla v||_2^2+\frac{\alpha}{2} \int_{\mathbb{R}^{3}}\left(u^2+v^2\right)\left(\left | u \right |+\left | v \right | \right){\rm d}x\\
 &\ \ \ -\frac{\alpha }{16\pi}\int_{\mathbb{R}^{3}}\left(\nabla (\left | u \right |+\left | v \right | )^{2}\right) {\rm d}x-\frac{1}{2}\left(||u||^{3}_{3}+||v||_{3}^{3}\right)-\beta\int_{\mathbb{R}^{3}}|u|^{\frac{3}{2}}|v|^{\frac{3}{2}}{\rm d}x\\
 &\ge \left(1-\frac{\alpha }{8\pi} \right)\left(||\nabla u||_2^2+||\nabla v||_2^2\right)+\frac{\alpha -1}{2}\left(||u||^{3}_{3}+||v||_{3}^{3}\right)+\left(\alpha-\beta\right) \int_{\mathbb{R}^{3}}|u|^{\frac{3}{2}}|v|^{\frac{3}{2}}{\rm d}x.
\end{align*}
Since $1\le \alpha<8\pi$, and $0<\beta<\alpha$, we can get $Q\left(u,v\right)>0$. At this point, the proof is complete.
\s{Proof of Theorem \ref{stability}}
In this section, we prove Theorem \ref{stability} following the classical arguments of \cite{Cazenave-6,Grillakis-14}. \\
$ Proof\ of \ Theorem$ \ref{stability}.
First of all,  we notice explicitly that  $G\left(a_{1},a_{2}\right)$ is invariant by translation, i.e. if $ \left ( u,v \right )\in G\left ( a_{1} ,a_{2}\right ) $ then also  $ \left ( u(\cdot-y) ,v(\cdot-y) \right )\in G\left ( a_{1} ,a_{2}\right ) $  for any $y \in \mathbb{R}^{3}$.
We argue by contradiction, assuming that there exist  $a_{1}$ and $a_{2}>0$ such that $G\left(a_{1},a_{2}\right)$ is not orbitally stable. This means that there is a $\epsilon_{0} >0$, and a sequence of initial $ \left(\Psi _{1}^{n} \left ( 0 \right ) ,\Psi _{2}^{n} \left ( 0 \right )\right)\subset  H^{1}_{r} \left(\mathbb{R}^{3}  \right)\times H^{1}_{r} \left(\mathbb{R}^{3}  \right)$ and $\left \{ t_{n}  \right \} \subset \mathbb{R^{+} } $ such that
$$
\inf _{\left ( u,v\right )\in G\left ( a_{1} ,a_{2}\right )  }\left \| (\Psi _{1}^{n} \left ( 0 \right ) ,\Psi _{2}^{n} \left ( 0 \right ))-\left ( u,v \right ) \right \|_{H^{1}}\longrightarrow 0
$$
and
\be\lab{contradiction}
\inf _{\left ( u,v\right )\in G\left ( a_{1} ,a_{2}\right )  }\left \| \left(\Psi _{1}^{n} \left ( \cdot,t_{n}  \right ) ,\Psi _{2}^{n} \left ( \cdot,t_{n}  \right )\right)-\left (u,v\right ) \right \|_{H^{1}} \ge \varepsilon_{0} .
\ee
Since by the conservation laws, the energy and the charge associated with $\Psi _{i}\left(,t\right)$ $i=1,2$. satisfies
$I\left(\Psi _{1}^{n}(\cdot,t_{n}),\Psi _{2}^{n}(\cdot,t_{n})\right)=I\left(\Psi _{1}^{n}\left(\cdot,0\right),\Psi _{2}^{n}\left(\cdot,0\right)\right)$, and $\left \| \Psi _{i}^{n}(\cdot,t_{n}) \right \| _{2}^{2}=\left \| \Psi _{i}^{n}(0) \right \| _{2}^{2} $,  for $i=1,2$.
Define
$$
\tilde{\Psi _{i}^{n}}(\cdot,t_{n})=\frac{\Psi _{i}^{n}(\cdot,t_{n})}{\left \| \Psi _{i}^{n}(\cdot,t_{n}) \right \|_{2}^{2}  }a_{i}^{\frac{1}{2} }, \text{for} \quad i=1,2,
$$
we have
$$
\left \| \Psi _{i}^{n}(\cdot,t_{n}) \right \| _{2}^{2}=a_{i},  \text{for} \quad i=1,2, \text{and} \quad I\left ( \tilde{\Psi _{1}^{n}}, \tilde{\Psi _{2}^{n}}  \right ) =m\left ( a_{1},a_{2}   \right ) + o _{n}\left ( 1 \right ) .
$$
So we find a minimizing sequence $( \tilde{\Psi _{1}^{n}}, \tilde{\Psi _{2}^{n}}$  with respect to $m\left ( a_{1},a_{2}  \right )$. However, according to theorem \ref{mainresult} (1)(2), the minimizing sequence is precompact (up to translation) in $H_{r}^{1}(\mathbb{R}^{3}) \times H_{r}^{1}(\mathbb{R}^{3})$,  which contradicts (\ref{contradiction}). The proof is completed.

\end{document}